\documentclass[12pt]{article}
\usepackage{amssymb}

\begin{document}

\baselineskip 18pt

\title{On two-parameter deformations of $osp(1|2)^{(1)}$}
\author{Liu Zhao${}^{1,2}$\thanks{%
Royal Society visiting fellow} \hspace{1cm} Xiang-Mao Ding${}^{3,4}$ \\
${}^1$ Institute of Modern Physics, Northwest University, Xian 710069, China%
\\
${}^2$ Department of Mathematics, University of York, York YO10 5DD, UK\\
${}^3$ Institute of Applied Mathematics,\\
Academy of Mathematics and Systems Science,\\
Academia Sinica, Beijing 100080,China\\
${}^4$ Institute of Theoretical Physics,\\
Academy of China, Beijing 100080, China}
\maketitle


\begin{center}
ABSTRACT

\begin{minipage}{5in}
An elliptic two-parameter deformation of the (universal enveloping superalgebra of)
affine Lie superalgebra $osp(1|2)^{(1)}$ is proposed in terms of free boson realization.
This deformed superalgebra is shown to fit in the framework of infinite Hopf family of
superalgebras, a generalization of the infinite Hopf family of algebras proposed
earlier by the authors. The trigonometric and rational degenerations are briefly discussed.
\end{minipage}
\end{center}

\vspace{1cm}

It has been becoming more and more evident that quantum affine \cite{D1,D2,D3}
and quantum Virasoro algebras \cite{FF} play some essential
roles in the theories of massive
integrable quantum fields in (1+1)-dimensions and in 2d statistical
mechanical systems off the criticality. For some systems, similar roles are
played by some yet more complicated algebraic structures, e.g. some
two-parameter deformations of affine Lie algebras. It has been realized that
there are several kinds of two-parameter deformations of affine Lie
algebras, including the standard elliptic quantum groups proposed by Felder
et al. \cite{felder1,felder2,felder3} and some different variants thereof
\cite{Foda1,Foda2,HY,f,hou2,jimbot,Konno,YZ1,Z1}. All these two-parameter algebras
fall into one of the following two classes: one is the quasi-triangular
quasi-Hopf algebra \cite{D4}--certain twists \cite{jimbot} of the standard
Hopf algebra structure \cite{D1,D2}--and the other is the so-called
infinite Hopf family of algebras \cite{f,hou2}. To us the relation
between the structures of quasi-triangular quasi-Hopf algebras and
infinite Hopf family of algebras still remains a mystery because
these two structures are defined respectively for algebras given
in different realizations: the quasi-triangular quasi-Hopf algebra
is introduced in the context of Yang-Baxter realization 
(or RS-realization \cite{RS})
and the co-structure is a deformation of the standard co-algebraic
structure of quantum affine algebras, while
the infinite Hopf family of algebras is introduced for algebras
given in the Drinfeld new current realization \cite{D3} only,
and the co-structure is a deformation of Drinfeld's new
co-structure for quantum affine algebras \cite{f,hou2}.

Among the above mentioned two-parameter deformations of affine algebras, we
are particularly interested in the algebras $\mathcal{E}_{q,p}(\widehat{g})$
studied in \cite{f,hou2} because, at level $c=1$, such algebras are
intimately related to the elliptic quantum $W$-algebras proposed by
B. Feigin and E.Frenkel \cite{FF}: the former are basically the algebra
of screening currents of the latter, with the introduction of some
auxiliary semi-simple currents and perhaps some dynamical shift.

In this article, we are aiming at proposing a generalization of
the algebra $\mathcal{E}_{q,p}(\widehat{g})$ to the
case of Lie superalgebra, with $g$ identified as the simplest
Lie superalgebra $osp(1|2)$. One of the reason to do such a work
is due to the fact that, although the elliptic quantum $W$%
-algebras have been fairly well understood \cite {FF},
their super-correspondences has
not been properly studied, even the structure of the simplest super Virasoro
algebra(if any) is still unknown. From the experiences of purely bosonic
case, it may be reasonable to expect that, the elliptic quantum super
Virasoro algebra, if exists, will possess an algebra of screening currents
which is based on the underlying affine Lie superalgebra $osp(1|2)^{(1)}$
and has a similar structure as that of $\mathcal{E}_{q,p}(\widehat{g}).$
Therefore, the study of a two-parameter elliptic quantum superalgebra might
hopefully give some hints about the unknown elliptic quantum super Virasoro
algebra. On the other hand, the study of such an algebra is of its own
interests also: it will provide examples of super analogues of two-parameter
quantum deformed affine algebras which naturally fit in the framework of
vertex operator algebras studied enthusiastically by pure mathematicians.
The existence, structure and representation theories of such algebras are
all worth studying.

\vspace{1pt}

The logic we shall be going along with is as follows. First we show that we
can close a set of currents into an associative superalgebra starting from
some free boson expressions. Then we show that this current superalgebra is
a particular case (level $c=1$) of some more general, abstract superalgebra $%
\mathcal{E}_{q,p}(osp(1|2)^{(1)})$ associated with $osp(1|2)^{(1)}$. Next,
we prove that the superalgebra $\mathcal{E}_{q,p}(osp(1|2)^{(1)})$ has a
well-defined co-structure, which can be regarded as a super analogue of the
structure of the so-called infinite Hopf family of algebras defined earlier by
the authors for $\mathcal{E}_{q,p}(\widehat{g})$. Last, we show that the
superalgebra $\mathcal{E}_{q,p}(osp(1|2)^{(1)})$ has some interesting
scaling limit, denoted as $\mathcal{A}_{\hbar ,\eta }(osp(1|2)^{(1)})$,
which is the analogue of $\mathcal{A}_{\hbar ,\eta }(\widehat{g}),$ the
trigonometric degeneration of $\mathcal{E}_{q,p}(\widehat{g})$.

\vspace{1pt}

To begin with, we introduce the following Heisenberg algebra $\mathcal{H}$
with generators \{$a_{i},$ $i\in Z$\} and relations

\begin{eqnarray*}
\lbrack a_{n},a_{m}\rbrack &=&\frac{1}{n}\left( q^{n}-q^{-n}\right) \left(
\left( qp\right) ^{n}-\left( qp\right) ^{-n}\right) \left(
p^{n}+p^{-n}-1\right) \delta _{n+m,0},\qquad (n\neq 0) \\
\lbrack P,Q\rbrack &=&1.
\end{eqnarray*}

\vspace{1pt}

\noindent Let

\vspace{1pt}

\[
s_{n}^{+}=\frac{a_{n}}{q^{n}-q^{-n}},\qquad s_{n}^{-}=\frac{a_{n}}{\left(
qp\right) ^{n}-\left( qp\right) ^{-n}}.
\]

\vspace{1pt}

\noindent Define

\vspace{1pt}

\begin{eqnarray}
\varphi (z)=\sum_{n\neq 0}s_{n}^{+}z^{-n},\qquad \psi (z)=\sum_{n\neq
0}s_{n}^{-}z^{-n},  \label{a00}
\end{eqnarray}

\vspace{1pt}

\noindent we have

\vspace{1pt}
\begin{eqnarray*}
\langle \varphi (z)\varphi (w)\rangle &=&-\sum_{n=0}^{\infty }\frac{1}{n}%
\frac{\left( \left( qp\right) ^{n}-\left( qp\right) ^{-n}\right) \left(
p^{n}+p^{-n}-1\right) }{q^{n}-q^{-n}}\left( \frac{w}{z}\right) ^{n}, \\
\langle \psi (z)\psi (w)\rangle &=&-\sum_{n=0}^{\infty }\frac{1}{n}\frac{%
\left( q^{n}-q^{-n}\right) \left( p^{n}+p^{-n}-1\right) }{\left( qp\right)
^{n}-\left( qp\right) ^{-n}}\left( \frac{w}{z}\right) ^{n}, \\
\langle \varphi (z)\psi (w)\rangle &=&-\sum_{n=0}^{\infty }\frac{1}{n}\left(
p^{n}+p^{-n}-1\right) \left( \frac{w}{z}\right) ^{n}.
\end{eqnarray*}

\vspace{1pt}

\noindent Therefore,

\vspace{1pt}

\begin{eqnarray*}
\exp \left\{ \langle \varphi (z)\varphi (w)\rangle \right\} &=&\exp \left\{
-\sum_{n=1}^{\infty }\frac{1}{n}\frac{\left( \left( qp\right) ^{n}-\left(
qp\right) ^{-n}\right) \left( p^{n}+p^{-n}-1\right) }{q^{n}-q^{-n}}\left(
\frac{w}{z}\right) ^{n}\right\} \\
&=&\frac{\left( \frac{w}{z}p^{-2}|q^{2}\right) _{\infty }\left( \frac{w}{z}%
q^{2}p|q^{2}\right) _{\infty }\left( \frac{w}{z}|q^{2}\right) _{\infty }}{%
\left( \frac{w}{z}(qp)^{2}|q^{2}\right) _{\infty }\left( \frac{w}{z}%
p^{-1}|q^{2}\right) _{\infty }\left( \frac{w}{z}q^{2}|q^{2}\right) _{\infty }%
};
\end{eqnarray*}

\vspace{1pt}

\begin{eqnarray*}
\exp \left\{ \langle \psi (z)\psi (w)\rangle \right\} &=&\exp \left\{
-\sum_{n=1}^{\infty }\frac{1}{n}\frac{\left( q^{n}-q^{-n}\right) \left(
p^{n}+p^{-n}-1\right) }{\left( qp\right) ^{n}-\left( qp\right) ^{-n}}\left(
\frac{w}{z}\right) ^{n}\right\} \\
&=&\frac{\left( \frac{w}{z}p^{2}|\left( qp\right) ^{2}\right) _{\infty
}\left( \frac{w}{z}\left( qp\right) ^{2}p^{-1}|\left( qp\right) ^{2}\right)
_{\infty }\left( \frac{w}{z}|\left( qp\right) ^{2}\right) _{\infty }}{\left(
\frac{w}{z}\left( qp\right) ^{2}p^{-2}|\left( qp\right) ^{2}\right) _{\infty
}\left( \frac{w}{z}p|\left( qp\right) ^{2}\right) _{\infty }\left( \frac{w}{z%
}\left( qp\right) ^{2}|\left( qp\right) ^{2}\right) _{\infty }};
\end{eqnarray*}

\vspace{1pt}

\begin{eqnarray*}
\exp \left\{ \langle \varphi (z)\psi (w)\rangle \right\} &=&\exp \left\{
-\sum_{n=1}^{\infty }\frac{1}{n}\left( p^{n}+p^{-n}-1\right) \left( \frac{w}{%
z}\right) ^{n}\right\} \\
&=&\exp \left\{ \log \left( 1-\frac{w}{z}p\right) +\log \left( 1-\frac{w}{z}%
p^{-1}\right) -\log \left( 1-\frac{w}{z}\right) \right\} \\
&=&\frac{\left( 1-\frac{w}{z}p\right) \left( 1-\frac{w}{z}p^{-1}\right) }{%
\left( 1-\frac{w}{z}\right) },
\end{eqnarray*}

\noindent where $(z|q)_\infty = \prod_{n=0}^{\infty}(1-zq^n)$. Now define

\vspace{1pt}

\begin{eqnarray}
S^{+}(z)=:\exp \left[ \varphi (z)\right] :,\qquad S^{-}(z)=:\exp \left[
-\psi (z)\right] :,  \label{a01}
\end{eqnarray}

\vspace{1pt}

\noindent where : \hspace{0.1cm} : means the usual normal ordering of
bosonic oscillators, we have

\vspace{1pt}
\begin{eqnarray*}
S^{+}(z)S^{+}(w) &=&\exp \left\{ \langle \varphi (z)\varphi (w)\rangle
\right\} :S^{+}(z)S^{+}(w):, \\
S^{-}(z)S^{-}(w) &=&\exp \left\{ \langle \psi (z)\psi (w)\rangle \right\}
:S^{-}(z)S^{-}(w):, \\
S^{+}(z)S^{-}(w) &=&\exp \left\{ -\langle \varphi (z)\psi (w)\rangle
\right\} :S^{+}(z)S^{-}(w):.
\end{eqnarray*}

\noindent Notice that the above bosonic expressions do not contain the
contribution of the zero mode and therefore do not live in the complete Fock
space corresponding to the Heisenberg algebra. Adding the zero mode
contributions, we now introduce

\begin{eqnarray}
E(z)=e^{Q}z^{P}S^{+}(z),\qquad F(z)=e^{-Q}z^{-P}S^{-}(z).  \label{a02}
\end{eqnarray}

\noindent Then, it is easy to calculate that

\begin{eqnarray*}
E(z)E(w) &=&\left( z-w\right) \frac{\left( \frac{w}{z}p^{-2}|q^{2}\right)
_{\infty }\left( \frac{w}{z}q^{2}p|q^{2}\right) _{\infty }}{\left( \frac{w}{z%
}(qp)^{2}|q^{2}\right) _{\infty }\left( \frac{w}{z}p^{-1}|q^{2}\right)
_{\infty }}:E(z)E(w):, \\
F(z)F(w) &=&\left( z-w\right) \frac{\left( \frac{w}{z}p^{2}|\left( qp\right)
^{2}\right) _{\infty }\left( \frac{w}{z}\left( qp\right) ^{2}p^{-1}|\left(
qp\right) ^{2}\right) _{\infty }}{\left( \frac{w}{z}\left( qp\right)
^{2}p^{-2}|\left( qp\right) ^{2}\right) _{\infty }\left( \frac{w}{z}p|\left(
qp\right) ^{2}\right) _{\infty }}:F(z)F(w):, \\
E(z)F(w) &=&\frac{\left( z-w\right) }{z^{2}\left( 1-\frac{w}{z}p\right)
\left( 1-\frac{w}{z}p^{-1}\right) }:E(z)F(w):.
\end{eqnarray*}

\noindent For $z=w$, the OPEs $E(z)E(w)$ and $F(z)F(w)$ are both zero,
showing that the currents $E(z),F(w)$ are essentially fermionic. Turning the
above equations into commutator type, we have

\vspace{1pt}

\begin{eqnarray*}
E(z)E(w)+\frac{\theta _{q^{2}}\left( \frac{w}{z}p^{-2}\right) \theta
_{q^{2}}\left( \frac{z}{w}p^{-1}\right) }{\theta _{q^{2}}\left( \frac{z}{w}%
p^{-2}\right) \theta _{q^{2}}\left( \frac{w}{z}p^{-1}\right) }E(w)E(z) &=&0,
\\
F(z)F(w)+\frac{\theta _{(qp)^{2}}\left( \frac{w}{z}p^{2}\right) \theta
_{(qp)^{2}}\left( \frac{z}{w}p\right) }{\theta _{(qp)^{2}}\left( \frac{z}{w}%
p^{2}\right) \theta _{(qp)^{2}}\left( \frac{w}{z}p\right) }F(w)F(z) &=&0,
\end{eqnarray*}

\vspace{1pt}

\noindent or alternatively,

\vspace{1pt}

\begin{eqnarray*}
E(z)E(w)+p\frac{\theta _{q^{2}}\left( \frac{w}{z}p^{-2}\right) \theta
_{q^{2}}\left( \frac{w}{z}p\right) }{\theta _{q^{2}}\left( \frac{w}{z}%
p^{2}\right) \theta _{q^{2}}\left( \frac{w}{z}p^{-1}\right) }E(w)E(z) &=&0,
\\
F(z)F(w)+p^{-1}\frac{\theta _{(qp)^{2}}\left( \frac{w}{z}p^{2}\right) \theta
_{(qp)^{2}}\left( \frac{w}{z}p^{-1}\right) }{\theta _{(qp)^{2}}\left( \frac{w%
}{z}p^{-2}\right) \theta _{(qp)^{2}}\left( \frac{w}{z}p\right) }F(w)F(z)
&=&0,
\end{eqnarray*}

\vspace{1pt}

\noindent where $\theta _{q}(z)=(z|q)_{\infty }(z^{-1}q|q)_{\infty
}(q|q)_{\infty }$ is essentially the usual Jacobi $\theta$-function\footnote{
To be precise, we have $\theta _{1}(u,\tau)=-iq^{1/8}z^{-1/2}\theta _{q}(z)$,
where $z=e^{2i \pi u}, q=e^{2\pi i\tau }$.}. Meanwhile, we also have

\vspace{1pt}

\begin{eqnarray*}
\left\{ E(z),F(w)\right\} &=&\frac{z-w}{(p-p^{-1})zw}\left\{ \delta \left(
\frac{z}{wp}\right) -\delta \left( \frac{w}{zp}\right) \right\} :E(z)F(w): \\
&=&\frac{1}{(p-p^{-1})}\left\{ \frac{p-1}{wp}\delta \left( \frac{z}{wp}%
\right) +\frac{p-1}{zp}\delta \left( \frac{w}{zp}\right) \right\} :E(z)F(w):
\\
&=&\frac{1}{p^{1/2}+p^{-1/2}}\left\{ (wp^{1/2})^{-1}\delta \left( \frac{z}{wp%
}\right) +(zp^{1/2})^{-1}\delta \left( \frac{w}{zp}\right) \right\}
:E(z)F(w):.
\end{eqnarray*}

\vspace{1pt}

\noindent Therefore, defining

\vspace{1pt}

\[
H^{\pm }(z)=z^{-1}:E(zp^{\pm 1/2})F(zp^{\mp 1/2}):,
\]

\vspace{1pt}

\noindent we get

\[
\left\{ E(z),F(w)\right\} =\frac{1}{p^{1/2}+p^{-1/2}}\left\{ \delta \left(
\frac{z}{wp}\right) H^{+}(wp^{1/2})+\delta \left( \frac{w}{zp}\right)
H^{-}(zp^{1/2})\right\}
\]

\noindent and

\begin{eqnarray*}
H^{+}(z)E(w) &=&p\frac{\theta _{q^{2}}\left( \frac{w}{z}p^{-2}\cdot
p^{-1/2}\right) \theta _{q^{2}}\left( \frac{w}{z}p\cdot p^{-1/2}\right) }{%
\theta _{q^{2}}\left( \frac{w}{z}p^{2}\cdot p^{-1/2}\right) \theta
_{q^{2}}\left( \frac{w}{z}p^{-1}\cdot p^{-1/2}\right) }E(w)H^{+}(z), \\
H^{-}(z)E(w) &=&p\frac{\theta _{q^{2}}\left( \frac{w}{z}p^{-2}\cdot
p^{1/2}\right) \theta _{q^{2}}\left( \frac{w}{z}p\cdot p^{1/2}\right) }{%
\theta _{q^{2}}\left( \frac{w}{z}p^{2}\cdot p^{1/2}\right) \theta
_{q^{2}}\left( \frac{w}{z}p^{-1}\cdot p^{1/2}\right) }E(w)H^{+}(z), \\
H^{+}(z)F(w) &=&p^{-1}\frac{\theta _{(qp)^{2}}\left( \frac{w}{z}p^{2}\cdot
p^{-1/2}\right) \theta _{(qp)^{2}}\left( \frac{w}{z}p^{-1}\cdot
p^{-1/2}\right) }{\theta _{(qp)^{2}}\left( \frac{w}{z}p^{-2}\cdot
p^{-1/2}\right) \theta _{(qp)^{2}}\left( \frac{w}{z}p\cdot p^{-1/2}\right) }%
F(w)H^{+}(z), \\
H^{-}(z)F(w) &=&p^{-1}\frac{\theta _{(qp)^{2}}\left( \frac{w}{z}p^{2}\cdot
p^{1/2}\right) \theta _{(qp)^{2}}\left( \frac{w}{z}p^{-1}\cdot
p^{1/2}\right) }{\theta _{(qp)^{2}}\left( \frac{w}{z}p^{-2}\cdot
p^{1/2}\right) \theta _{(qp)^{2}}\left( \frac{w}{z}p\cdot p^{1/2}\right) }%
F(w)H^{-}(z), \\
H^{\pm }(z)H^{\pm }(w) &=&\frac{\theta _{q^{2}}\left( \frac{w}{z}%
p^{-2}\right) \theta _{q^{2}}\left( \frac{w}{z}p\right) }{\theta
_{q^{2}}\left( \frac{w}{z}p^{2}\right) \theta _{q^{2}}\left( \frac{w}{z}%
p^{-1}\right) }\frac{\theta _{(qp)^{2}}\left( \frac{w}{z}p^{2}\right) \theta
_{(qp)^{2}}\left( \frac{z}{w}p\right) }{\theta _{(qp)^{2}}\left( \frac{z}{w}%
p^{2}\right) \theta _{(qp)^{2}}\left( \frac{w}{z}p\right) }H^{\pm }(w)H^{\pm
}(z), \\
H^{+}(z)H^{-}(w) &=&\frac{\theta _{q^{2}}\left( \frac{w}{z}p^{-2}\cdot
p^{-1}\right) \theta _{q^{2}}\left( \frac{w}{z}p\cdot p^{-1}\right) }{\theta
_{q^{2}}\left( \frac{w}{z}p^{2}\cdot p^{-1}\right) \theta _{q^{2}}\left(
\frac{w}{z}p^{-1}\cdot p^{-1}\right) }\frac{\theta _{(qp)^{2}}\left( \frac{w%
}{z}p^{2}\cdot p\right) \theta _{(qp)^{2}}\left( \frac{z}{w}p\cdot p\right)
}{\theta _{(qp)^{2}}\left( \frac{z}{w}p^{2}\cdot p\right) \theta
_{(qp)^{2}}\left( \frac{w}{z}p\cdot p\right) }H^{-}(w)H^{+}(z).
\end{eqnarray*}

\noindent Notice that unlike $E(z)$ and $F(z)$, the currents $H^{\pm }(z)$
are bosonic. For later reference, we introduce the Grassmann parity operator
$\pi $ such that

\begin{eqnarray*}
\pi\lbrack E(z) \rbrack &=& \pi\lbrack F(z) \rbrack =1, \\
\pi\lbrack H^{\pm}(z) \rbrack &=& 0.
\end{eqnarray*}

\noindent Then we see that the currents $H^{\pm }(z),$ $E(z)$ and $F(z)$
close into a superalgebra with $Z_{2}$ gradation provided by the parity
operator $\pi $.

\noindent \textbf{Definition 1}: \emph{The superalgebra }$\mathcal{E}%
_{q,p}(osp(1|2)^{(1)})$\emph{\ is a }$Z_{2}$\emph{-graded associative
algebra generated by the unit }$1$\emph{, coefficients of the formal power
series }$H^{\pm }(z)$\emph{\ (invertible), }$E(z)$\emph{\ and }$F(z)$\emph{\
in }$z$\emph{\ and the central element }$c$\emph{\ with relations}

\begin{eqnarray*}
H^{+}(z)E(w) &=&p\frac{\theta _{q^{2}}\left( \frac{w}{z}p^{-2}\cdot
p^{-c/2}\right) \theta _{q^{2}}\left( \frac{w}{z}p\cdot p^{-c/2}\right) }{%
\theta _{q^{2}}\left( \frac{w}{z}p^{2}\cdot p^{-c/2}\right) \theta
_{q^{2}}\left( \frac{w}{z}p^{-1}\cdot p^{-c/2}\right) }E(w)H^{+}(z), \\
H^{-}(z)E(w) &=&p\frac{\theta _{q^{2}}\left( \frac{w}{z}p^{-2}\cdot
p^{c/2}\right) \theta _{q^{2}}\left( \frac{w}{z}p\cdot p^{c/2}\right) }{%
\theta _{q^{2}}\left( \frac{w}{z}p^{2}\cdot p^{c/2}\right) \theta
_{q^{2}}\left( \frac{w}{z}p^{-1}\cdot p^{c/2}\right) }E(w)H^{+}(z), \\
H^{+}(z)F(w) &=&p^{-1}\frac{\theta _{\tilde{q}^{2}}\left( \frac{w}{z}%
p^{2}\cdot p^{-c/2}\right) \theta _{\tilde{q}^{2}}\left( \frac{w}{z}%
p^{-1}\cdot p^{-c/2}\right) }{\theta _{\tilde{q}^{2}}\left( \frac{w}{z}%
p^{-2}\cdot p^{-c/2}\right) \theta _{\tilde{q}^{2}}\left( \frac{w}{z}p\cdot
p^{-c/2}\right) }F(w)H^{+}(z), \\
H^{-}(z)F(w) &=&p^{-1}\frac{\theta _{\tilde{q}^{2}}\left( \frac{w}{z}%
p^{2}\cdot p^{c/2}\right) \theta _{\tilde{q}^{2}}\left( \frac{w}{z}%
p^{-1}\cdot p^{c/2}\right) }{\theta _{\tilde{q}^{2}}\left( \frac{w}{z}%
p^{-2}\cdot p^{c/2}\right) \theta _{\tilde{q}^{2}}\left( \frac{w}{z}p\cdot
p^{c/2}\right) }F(w)H^{-}(z), \\
H^{\pm }(z)H^{\pm }(w) &=&\frac{\theta _{q^{2}}\left( \frac{w}{z}%
p^{-2}\right) \theta _{q^{2}}\left( \frac{w}{z}p\right) }{\theta
_{q^{2}}\left( \frac{w}{z}p^{2}\right) \theta _{q^{2}}\left( \frac{w}{z}%
p^{-1}\right) }\frac{\theta _{\tilde{q}^{2}}\left( \frac{w}{z}p^{2}\right)
\theta _{\tilde{q}^{2}}\left( \frac{z}{w}p\right) }{\theta _{\tilde{q}%
^{2}}\left( \frac{z}{w}p^{2}\right) \theta _{\tilde{q}^{2}}\left( \frac{w}{z}%
p\right) }H^{\pm }(w)H^{\pm }(z), \\
H^{+}(z)H^{-}(w) &=&\frac{\theta _{q^{2}}\left( \frac{w}{z}p^{-2}\cdot
p^{-c}\right) \theta _{q^{2}}\left( \frac{w}{z}p\cdot p^{-c}\right) }{\theta
_{q^{2}}\left( \frac{w}{z}p^{2}\cdot p^{-c}\right) \theta _{q^{2}}\left(
\frac{w}{z}p^{-1}\cdot p^{-c}\right) }\frac{\theta _{\tilde{q}^{2}}\left(
\frac{w}{z}p^{2}\cdot p^{c}\right) \theta _{\tilde{q}^{2}}\left( \frac{z}{w}%
p\cdot p^{c}\right) }{\theta _{\tilde{q}^{2}}\left( \frac{z}{w}p^{2}\cdot
p^{c}\right) \theta _{\tilde{q}^{2}}\left( \frac{w}{z}p\cdot p^{c}\right) }%
H^{-}(w)H^{+}(z), \\
E(z)E(w) &=&-p\frac{\theta _{q^{2}}\left( \frac{w}{z}p^{-2}\right) \theta
_{q^{2}}\left( \frac{w}{z}p\right) }{\theta _{q^{2}}\left( \frac{w}{z}%
p^{2}\right) \theta _{q^{2}}\left( \frac{w}{z}p^{-1}\right) }E(w)E(z), \\
F(z)F(w) &=&-p^{-1}\frac{\theta _{\tilde{q}^{2}}\left( \frac{w}{z}%
p^{2}\right) \theta _{\tilde{q}^{2}}\left( \frac{w}{z}p^{-1}\right) }{\theta
_{\tilde{q}^{2}}\left( \frac{w}{z}p^{-2}\right) \theta _{\tilde{q}%
^{2}}\left( \frac{w}{z}p\right) }F(w)F(z), \\
\left\{ E(z),F(w)\right\} &=&\frac{1}{p^{1/2}+p^{-1/2}}\left\{ \delta \left(
\frac{z}{wp^{c}}\right) H^{+}(wp^{c/2})+\delta \left( \frac{w}{zp^{c}}%
\right) H^{-}(zp^{c/2})\right\} ,
\end{eqnarray*}

\noindent \emph{where }$\tilde{q}=qp^{c}$.\hfill$\square $

\vspace{1pt}

Obviously, we have

\vspace{1pt}

\noindent \textbf{Proposition 1}:\emph{\ Equations (\ref{a00}, \ref{a01},
\ref{a02}) give a free boson realization for the superalgebra }$\mathcal{E}%
_{q,p}(osp(1|2)^{(1)})$\emph{\ at level }$c=1$\emph{.}\hfill $\square $

To understand the properties of the superalgebra $\mathcal{E}%
_{q,p}(osp(1|2)^{(1)})$ at levels other than $c=1$, one natural way is to
consider its co-structure. Since the superalgebra $\mathcal{E}%
_{q,p}(osp(1|2)^{(1)})$ has some similar properties as those of the $\mathcal{E}%
_{q,p}(\widehat{g})$, e.g. the Cartan involution is broken due to the
different periods of the structure functions for the $EE$ and $FF$
relations, we hope that the structure of infinite Hopf family of algebras
for the latter also holds for $\mathcal{E}_{q,p}(osp(1|2)^{(1)})$. It is
indeed so, however, the definition of infinite Hopf family of algebras has
to be modified into an infinite Hopf family of superalgebras, as is expected
naturally. This generalized co-structure for $\mathcal{E}%
_{q,p}(osp(1|2)^{(1)})$ will be given in the following proposition.

\vspace{1pt}

Let us first prepare some notations. Let $\mathcal{J}$ be an additive
semigroup which may be identified with the set of non-negative integer
numbers. $c_{n}\in \mathcal{J}$ are elements of $\mathcal{J}$. Let $%
q^{(0)}=q $ and define $q^{(n+1)}=q^{(n)}p^{c_{n}}$ iteratively. We set $%
\mathcal{A}_{n}=\mathcal{E}_{q^{(n)},p}(osp(1|2)^{(1)})$ whose generator are
denoted $H^{\pm }(z;q^{(n)}),$ $E(z;q^{(n)}),$ $F(z;q^{(n)})$ and $c_{n}$
respectively. The generating relations for each $\mathcal{A}_{n}$ are
nothing but those of $\mathcal{E}_{q,p}(osp(1|2)^{(1)})$ with parameters $q,%
\tilde{q}$ replaced by $q^{(n)},q^{(n+1)}$ respectively.

\vspace{1pt}

Let $\{v_{i}^{(n)},~i=1,~...,~\mathrm{\dim }(\mathcal{A}_{n})\}$ be a basis
of $\mathcal{A}_{n}$. The maps

\begin{eqnarray*}
\tau_n^{\pm}: \mathcal{A}_n &\rightarrow& \mathcal{A}_{n \pm 1} \\
v^{(n)}_i &\mapsto& v^{(n\pm 1)}_i
\end{eqnarray*}

\noindent are morphisms from $\mathcal{A}_n$ to $\mathcal{A}_{n \pm 1}$. For
any two integers $n,~m$ with $n<m$, we can specify a pair of morphisms

\begin{eqnarray*}
&&Mor(\mathcal{A}_{m},~\mathcal{A}_{n})\ni \tau ^{(m,n)}\equiv \tau
_{m-1}^{+}...\tau _{n+1}^{+}\tau _{n}^{+}:~~\mathcal{A}_{n}\rightarrow
\mathcal{A}_{m}, \\
&&Mor(\mathcal{A}_{n},~\mathcal{A}_{m})\ni \tau ^{(n,m)}\equiv \tau
_{n+1}^{-}...\tau _{m-1}^{-}\tau _{m}^{-}:~~\mathcal{A}_{m}\rightarrow
\mathcal{A}_{n}
\end{eqnarray*}

\noindent with $\tau ^{(m,n)}\tau ^{(n,m)}=id_{m},~\tau ^{(n,m)}\tau
^{(m,n)}=id_{n}$. Clearly the morphisms $\tau ^{(m,n)},~n,m\in Z$ satisfy
the associativity condition $\tau ^{(m,p)}\tau ^{(p,n)}=\tau ^{(m,n)}$ and
thus make the family of superalgebras $\{\mathcal{A}_{n},~n\in Z\}$ into a
category.

\vspace{1pt}

The following definition is a straightforward generalization of the
structure of infinite Hopf family of algebras originally presented in
\cite{f,hou2}:

\vspace{1pt}

\noindent \textbf{Definition 2}: \emph{The category of superalgebras }$\{%
\mathcal{A}_{n},~\{\tau ^{(n,m)}\},~n,m\in Z\}$\emph{\ is called an infinite
Hopf family of superalgebras if on each object }$\mathcal{A}_{n}$\emph{\ of
the category one can define the morphisms }$\Delta _{n}^{+}:\mathcal{A}%
_{n}\rightarrow \mathcal{A}_{n}\otimes \mathcal{A}_{n+1}$\emph{, }$\Delta
_{n}^{-}:\mathcal{A}_{n}\rightarrow \mathcal{A}_{n-1}\otimes \mathcal{A}_{n}$%
\emph{, }$\epsilon _{n}:\mathcal{A}_{n}\rightarrow C$\emph{\ and
antimorphisms }$S_{n}^{\pm }:\mathcal{A}_{n}\rightarrow \mathcal{A}_{n\pm 1}$%
\emph{\ such that the following axioms hold,}

\begin{itemize}
\item  $(\epsilon _{n}\otimes id_{n+1})\circ \Delta _{n}^{+}=\tau
_{n}^{+},~(id_{n-1}\otimes \epsilon _{n})\circ \Delta _{n}^{-}=\tau _{n}^{-}$%
\emph{\ \hfill (a1)}

\item  $m_{n+1}\circ (S_{n}^{+}\otimes id_{n+1})\circ \Delta
_{n}^{+}=\epsilon _{n+1}\circ \tau _{n}^{+},~m_{n-1}\circ (id_{n-1}\otimes
S_{n}^{-})\circ \Delta _{n}^{-}=\epsilon _{n-1}\circ \tau _{n}^{-}$\emph{\
\hfill (a2)}

\item  $(\Delta _{n}^{-}\otimes id_{n+1})\circ \Delta
_{n}^{+}=(id_{n-1}\otimes \Delta _{n}^{+})\circ \Delta _{n}^{-}$\emph{\
\hfill (a3)}
\end{itemize}

\emph{\noindent in which }$m_{n}$\emph{\ is the (super)multiplication for }$%
\mathcal{A}_{n}$\emph{.} \hfill$\square $

\noindent Notice that throughout this article, the symbol $\otimes $ denotes
a graded direct product, or \emph{direct super-product}, obeying, e.g. for
elements $A,B,C,D$ with definite Grassmann parity,

\begin{eqnarray*}
(A \otimes B)(C \otimes D) = (-1)^{\pi(B)\pi(C)} AC \otimes BD.
\end{eqnarray*}

\emph{\vspace{1pt}}

\noindent \textbf{Proposition 2}: \emph{The family of superalgebras }$\{%
\mathcal{A}_{n},$\emph{\ }$n\in Z\}$\emph{\ form an Infinite Hopf family of
algebras with comultiplications }$\Delta _{n}^{\pm },$\emph{\ counits }$%
\epsilon _{n}$\emph{\ and antipodes }$S_{n}^{\pm }$\emph{\ given as follows,}

\begin{itemize}
\item  \emph{\vspace{1pt}the comultiplications }$\Delta _{n}^{\pm }$\emph{: }
\begin{eqnarray*}
\Delta _{n}^{+}c_{n} &=&c_{n}+c_{n+1}, \\
\Delta _{n}^{+}H^{+}(z;q^{(n)}) &=&H^{+}(zp^{c_{n+1}/2};q^{(n)})\otimes
H^{+}(zp^{-c_{n}/2};q^{(n+1)}), \\
\Delta _{n}^{+}H^{-}(z;q^{(n)}) &=&-H^{-}(zp^{-c_{n+1}/2};q^{(n)})\otimes
H^{-}(zp^{c_{n}/2};q^{(n+1)}), \\
\Delta _{n}^{+}E(z;q^{(n)}) &=&E(z;q^{(n)})\otimes
1-H^{-}(zp^{c_{n}/2};q^{(n)})\otimes E(zp^{c_{n}};q^{(n+1)}), \\
\Delta _{n}^{+}F(z;q^{(n)}) &=&1\otimes
F(z;q^{(n+1)})+F(zp^{c_{n+1}};q^{(n)})\otimes
H^{+}(zp^{c_{n+1}/2};q^{(n+1)}), \\
&& \\
\Delta _{n}^{-}c_{n} &=&c_{n-1}+c_{n}, \\
\Delta _{n}^{-}H^{+}(z;q^{(n)}) &=&H^{+}(zp^{c_{n}/2};q^{(n-1)})\otimes
H^{+}(zp^{-c_{n-1}/2};q^{(n)}), \\
\Delta _{n}^{-}H^{-}(z;q^{(n)}) &=&-H^{-}(zp^{-c_{n}/2};q^{(n-1)})\otimes
H^{-}(zp^{c_{n-1}/2};q^{(n)}), \\
\Delta _{n}^{-}E(z;q^{(n)}) &=&E(z;q^{(n-1)})\otimes
1-H^{-}(zp^{c_{n-1}/2};q^{(n-1)})\otimes E(zp^{c_{n-1}};q^{(n)}), \\
\Delta _{n}^{-}F(z;q^{(n)}) &=&1\otimes
F(z;q^{(n)})+F(zp^{c_{n}};q^{(n-1)})\otimes H^{+}(zp^{c_{n}/2};q^{(n)});
\end{eqnarray*}

\item  \emph{the counits }$\epsilon _{n}$\emph{: }
\begin{eqnarray*}
\epsilon _{n}(c_{n}) &=&0, \\
\epsilon _{n}(1_{n}) &=&1, \\
\epsilon _{n}(H^{\pm }(z;q^{(n)})) &=&1, \\
\epsilon _{n}(E(z;q^{(n)})) &=&0, \\
\epsilon _{n}(F(z;q^{(n)})) &=&0;
\end{eqnarray*}

\item  \emph{the antipodes }$S_{n}^{\pm }$\emph{:}

\begin{eqnarray*}
S_{n}^{\pm }c_{n} &=&-c_{n\pm 1}, \\
S_{n}^{\pm }H^{+}(z;q^{(n)}) &=&\lbrack H^{+}(z;q^{(n\pm 1)})\rbrack ^{-1},
\\
S_{n}^{\pm }H^{-}(z;q^{(n)}) &=&\lbrack H^{-}(z;q^{(n\pm 1)})\rbrack ^{-1},
\\
S_{n}^{\pm }E(z;q^{(n)}) &=&-H^{-}(zp^{-c_{n\pm 1}/2};q^{(n\pm
1)})^{-1}E(zp^{-c_{n\pm 1}};q^{(n\pm 1)}), \\
S_{n}^{\pm }F(z;q^{(n)}) &=&F(zp^{-c_{n\pm 1}};q^{(n\pm
1)})H^{+}(zp^{-c_{n\pm 1}/2};q^{(n\pm 1)})^{-1}.
\end{eqnarray*}
\end{itemize}

\hfill$\square $

Let us stress that, among the defining relations of the superalgebra $%
\mathcal{E}_{q,p}(osp(1|2)^{(1)})$, the unusual signature in between the two
$\delta$-function terms in the relation containing the
anti-commutator of $E(z)$ and $F(w)$ is superficial: we can always
replace $H^-(z)$ by $-H^-(z)$ and e.g. $F(z)$ by $F(z)(p^{1/2}+p^{-1/2})/(p-p^{-1})$
and everything looks standard as in the usual q-affine algebra
case.

It is remarkable that the comultiplication $\Delta _{n}^{+}$ can
be applied iteratively onto $\mathcal{A}_{n},$ so that beginning
from the $c=1$ realization one can obtain a realization of higher
$c\in Z_{+}.$

If we re-parameterize the parameters $q,$ $p$ and $z$ as

\[
q=e^{\epsilon /\eta },p=e^{\epsilon \hbar },z=e^{i\epsilon u}
\]

\noindent and taking the scaling limit $\epsilon \rightarrow 0,$ the
superalgebra $\mathcal{E}_{q,p}(osp(1|2)^{(1)})$ degenerates into

\vspace{1pt}

\begin{eqnarray*}
H^{\pm }(u)E(v) &=&\frac{\sin 2\pi \eta (u-v+2\hbar \pm \hbar c/2)\sin 2\pi
\eta (u-v-\hbar \pm \hbar c/2)}{\sin \pi \eta (u-v-2\hbar \pm \hbar c/2)\sin
2\pi \eta (u-v+\hbar \pm \hbar c/2)}E(v)H^{+}(u), \\
H^{\pm }(u)F(v) &=&\frac{\sin 2\pi \eta ^{\prime }(u-v+2\hbar \mp \hbar
c/2)\sin 2\pi \eta ^{\prime }(u-v-\hbar \mp \hbar c/2)}{\sin 2\pi \eta
^{\prime }(u-v-2\hbar \mp \hbar c/2)\sin 2\pi \eta ^{\prime }(u-v+\hbar \mp
\hbar c/2)}F(v)H^{\pm }(u), \\
H^{\pm }(u)H^{\pm }(v) &=&\frac{\sin 2\pi \eta (u-v+2\hbar )\sin 2\pi \eta
(u-v-\hbar )}{\sin 2\pi \eta (u-v-2\hbar )\sin 2\pi \eta (u-v+\hbar )} \\
&&\times \frac{\sin 2\pi \eta ^{\prime }(u-v-2\hbar )\sin 2\pi \eta ^{\prime
}(u-v+\hbar )}{\sin 2\pi \eta ^{\prime }(u-v+2\hbar )\sin 2\pi \eta ^{\prime
}(u-v-\hbar )}H^{\pm }(v)H^{\pm }(u), \\
H^{+}(u)H^{-}(v) &=&\frac{\sin 2\pi \eta (u-v+2\hbar +\hbar c)\sin 2\pi \eta
(u-v-\hbar +\hbar c)}{\sin 2\pi \eta (u-v-2\hbar +\hbar c)\sin 2\pi \eta
(u-v+\hbar +\hbar c)} \\
&&\times \frac{\sin 2\pi \eta ^{\prime }(u-v-2\hbar -\hbar c)\sin 2\pi \eta
^{\prime }(u-v+\hbar -\hbar c)}{\sin 2\pi \eta ^{\prime }(u-v+2\hbar -\hbar
c)\sin 2\pi \eta ^{\prime }(u-v-\hbar -\hbar c)}H^{-}(v)H^{+}(u), \\
E(z)E(w) &=&-\frac{\sin 2\pi \eta (u-v+2\hbar )\sin 2\pi \eta (u-v-\hbar )}{%
\sin 2\pi \eta (u-v-2\hbar )\sin 2\pi \eta (u-v+\hbar )}E(w)E(z), \\
F(z)F(w) &=&-\frac{\sin 2\pi \eta ^{\prime }(u-v+2\hbar )\sin 2\pi \eta
^{\prime }(u-v-\hbar )}{\sin 2\pi \eta ^{\prime }(u-v-2\hbar )\sin 2\pi \eta
^{\prime }(u-v+\hbar )}F(w)F(z), \\
\left\{ E(u),F(v)\right\} &=&\frac{1}{2\hbar }\left\{ \delta \left(
u-v-\hbar c\right) H^{+}(v+\hbar c/2)+\delta \left( u-v+\hbar c\right)
H^{-}(u+\hbar c/2)\right\} ,
\end{eqnarray*}

\noindent where

\vspace{1pt}

\[
\frac{1}{\eta ^{\prime }}-\frac{1}{\eta }=\hbar c.
\]

\noindent This superalgebra is clearly an $osp(1|2)^{(1)}$ analogue of the
earlier studied algebras $\mathcal{A}_{\hbar ,\eta }(\widehat{g})$ and hence
we call it $\mathcal{A}_{\hbar ,\eta }(osp(1|2)^{(1)}).$ In the particular
case of $\eta \rightarrow 0$ this superalgebra further degenerates into the
super Yangian double $DY_{\hbar }(osp(1|2)^{(1)})$ -- the relations of which
(first introduced in \cite{ZYZ})
is just those of $\mathcal{A}_{\hbar ,\eta }(osp(1|2)^{(1)})$ but with all
the $\sin 2\pi \eta $ and $\sin 2\pi \eta ^{\prime }$ removed -- however
with a sign difference appeared in the last relation. The reason for this
sign difference has already been mentioned earlier in the context. This last
degeneration fully clarifies the connection of our superalgebras $\mathcal{E}%
_{q,p}(osp(1|2)^{(1)})$ and $\mathcal{A}_{\hbar ,\eta }(osp(1|2)^{(1)})$
with the underlying affine superalgebra $osp(1|2)^{(1)}$.

\vspace{1pt}

In closing, let us point out some related unsolved problems. The existence
of two parameter deformation of affine Lie (super)algebras with the
structure of infinite Hopf family of (super)algebras seems to be a universal
phenomenon, which means that there should be such a (super)algebra
associated with each underlying affine Lie (super)algebra. However, what we
have known about these (super)algebras is only a tiny top of an iceberg. We
know only a little about the structure theory and the representation theory.
The definition of these algebras themselves were only known for untwisted
affine Lie algebras associated with simply-laced Lie algebras and the
present article add to this picture the simplest affine Lie superalgebra $%
osp(1|2)^{(1)}$. It seems that there remains a lot of pure algebraic works
to do toward these (super)algebras.

\vspace{1pt}

On the other hand, physicists and/or applied mathematicians may be
particularly interested in the application aspect of these (super)algebras.
>From this point of view, we would like to mention the following problems
which we would like to see a solution:

\begin{itemize}
\item  Realizations of these (super)algebras other than the current
realization.

It is well known that for $q$-affine algebras and Yangian doubles
\cite{KT}, there are
mainly three different realizations which are proved to be connected
to each other: the current realization, Drinfeld realization in terms
of Laurent components of the currents and the Yang-Baxter realization
\cite{RS}. The last one is very important when physics applications
are considered because it relates the structure of quantum symmetry algebra and
the physical two-body $S$-matrix. So far we only know that, among the two parameter
deformed affine Lie (super)algebras with the structure of infinite Hipf family of
(super)algebras, only $\mathcal{A}_{\hbar,\eta}(sl(2)^{(1)})$ have a
Yang-Baxter realization which contains a``dynamical operator''
(spectral-shifting operator or weight vector of the underlying
Lie algebra) \cite{HY}. Attempts in obtaining Yang-Baxter
realizations for all other current algebras of this kind have
not lead to any success.

\item  The potential applications of these (super)algebras in physics
problems. Besides being related to the algebra of screening currents of the quantum
Virasoro and $W$-algebras, it would be interesting to see whether there is any
physical model which bears any of these algebras as the underlying quantum symmetry.
However, besides the case of $\mathcal{A}_{\hbar ,\eta}(sl(2)^{(1)})$, nobody has ever
been able to say a word on this possibility.

\item  Reconstruction of deformed $W$-algebras from the two parameter
deformed affine algebras. To us, this seems the most plausible route to seek
for physics applications, because these algebras are closely related to the
algebra of screening currents of the deformed $W$-algebras, and it is indeed
possible to reconstruct the $W$-algebras out of the screening currents. This
problem is particularly interesting when the superalgebra associated with $%
osp(1|2)^{(1)}$ is considered because we then will be able to gain some
knowledge about the deformed super Virasoro algebra -- an object expected
both for mathematical completeness and for physics applications!
\end{itemize}

\vspace{0.5cm}

\noindent \textbf{Acknowledgement}: L.Zhao would like to thank Niall MacKay
for hospitality at Dept. Appl. Math., Sheffield University and at Dept.
Math., Univ. of York during the preparation of this manuscript. The content
of this article has been presented at Sheffield University in an informal
seminar. This work is supported in part by the National Natural Science
Foundation of China.


\begin{thebibliography}{99}
\bibitem{D1}  Drinfeld, V.G., Hopf algebras and quantum Yang-Baxter
equations, \textit{Dokl. Akad. Nauk. SSSR 283 (1985) 1060}.

\bibitem{D2}  Drinfeld, V.G., Quantum groups, \textit{ICM Preceedings, 
New York, Berkeley (1986) 798}.

\bibitem{D3}  Drinfeld, V.G., New realizations of Yangian and quantum affine
algebras, \textit{Soviet. Math. Dokl. 26 (1988) 212}.

\bibitem{D4}  Drinfeld, V.G., Quasi-Hopf algebras, 
\textit{Liningrad Math. J. 1 (1990) 1419}.

\bibitem{FF}  Feigin B., Frenkel E., Quantum $W$-alebras and
Elliptiv algebras, \textit{q-alg/9508009}, \textit{Commun. Math. Phys.
178 (1996) 653}.

\bibitem{felder1}  Enriquez, B., Felder, G., Elliptic quantum groups $%
E_{\tau ,\eta }(sl_{2})$ and quasi-Hopf algebras, \textit{q-alg/9703018}.

\bibitem{felder2}  Felder, G., Conformal field theory and integrable systems
associated to elliptic curves, \textit{Proc. ICM Z\"{u}rich 1994, 1247,
Birkh\"{a}user (1994)}; Elliptic quantum groups, \textit{Proc. ICMP Paris
(1994) 2118, International Press (1995)}.

\bibitem{felder3}  Felder, G., Varchenko, A., On representation of
the elliptic quantum group $E_{\newline
\tau ,\eta }(sl_{2})$, \textit{q-alg/9601003}, \textit{Commun. Math. Phys.
181 (1996) 741}.

\bibitem{Foda1}  Foda, O., Iohara, K., Jimbo, M., Kedem, R., Miwa, T., Yan,
H., An elliptic quantum algebra for $\widehat{sl}_{2}$. \textit{Lett. Math.
Phys. 32 (1994) 259--268}.

\bibitem{Foda2}  Foda,O., Iohara,K., Jimbo,M., Kedem,R., Miwa,T., Yan, H.,
Notes on highest weight modules of the elliptic algebra $A_{p,q}(sl_{2})$.
\textit{Prog. Theoret. Phys., Supplement, 118 (1995) 1--34}.

\bibitem{ZYZ} Gould, M. and Zhang, Y.-Z., On super RS algebra and Drinfeld 
realization of quantum affine superalgebras, \textit{q-alg/9712011}, 
\textit{Lett. Math. Phys.44 (1998) 291}; Gould, M. and Zhang, Y.-Z.,
and Isaac,P.S., Casimir Invariants from Quasi-Hopf (Super)algebras, 
\textit{math.QA/9811062}, \textit{Commun. Math. Phys.}, in press.

\bibitem{HY}  Hou, B.-Y., Yang, W.-L., Dynamically twisted algebra $A_{q,p,%
\hat{\pi}}(\widehat{gl}_{2})$ as current algebra generalizing screening
currents of $q$-deformed Virasoro algebra, Preprint \textit{q-alg/9709024}.

\bibitem{f}  Hou, B.-Y., Zhao, L., Ding, X.-M., The algebra $\mathcal{A}%
_{\hbar ,\eta }(\hat{g})$ and infinite Hopf family of algebras, Preprint
\textit{q-alg/9703046}, \textit{J. Geom. Phys. 27 (1998) 249}.

\bibitem{hou2}  Hou, B.-Y., Zhao, L., Ding, X.-M., Infinite Hopf family of
elliptic algebras and bosonization, \textit{math/9801062}, 
\textit{J. Phys. A: Math. Gen. 32 (1999) 1951}.

\bibitem{jimbot}  Jimbo, M., Konno, H., Odake, S., Shiraishi, J., Quasi-Hopf
twistors for elliptic quantum groups, \textit{q-alg/9712029}.

\bibitem{KT}  Khoroshkin, S., Tolstoy, S., Yangian Double, 
\textit{Lett. Math. Phys.36 (1996) 373}.

\bibitem{Konno}  Konno, H., An elliptic algebra $U_{q,p}(\widehat{sl}_{2})$
and the fusion RSOS model, Preprint \textit{q-alg/9709013}.

\bibitem{RS}  Reshetikhin, N.Yu., Senemov-Tian-Shansky, M.A., Central
extensions of quantum current groups, \textit{Lett. Math. Phys. 19 
(1990) 133}.

\bibitem{YZ1} Zhang, Y.-Z. and Gould, M., Quasi-Hopf Superalgebras and 
Elliptic Quantm Supergroups, \textit{math.QA/9809156}, 
\textit{J. Math. Phys. 40 (1999) 5264}, and Zhang, Y.-Z., 
\textit{Prog. Theor. Phys. Suppl.135 (1999) 182}. 

\bibitem{Z1}  Zhao, L., B.-Y. Hou, Note on the algebra of screening currents
for the quantum deformed $W$ algebras, \textit{J. Phys. A: Math. Gen. 30
(1997) 7659}.

\end{thebibliography}
\end{document}